\documentclass[11pt]{amsart}
\usepackage{amsfonts}
\usepackage{amssymb}
\newtheorem{theorem}{{\bf Theorem}}

\newtheorem{corollary}[theorem]{{\bf Corollary}}

\newtheorem{problem}[theorem]{{\bf Problem}}
\newcommand{\pf}{{\bf Proof. \ }}
\newcommand{\epf}{{\square}}
\newcommand{\sone}{{\sf S}_1}

\newcommand{\sfin}{{\sf S}_{fin}}

\newcommand{\cpx}{{\sf C}_p(X)}
\newcommand{\ckx}{{\sf C}_k(X)}

\newcommand{\fplusx}{(2^X,{\sf F}^+)}
\newcommand{\zplusx}{(2^X,{\sf Z}^+)}

\newcommand{\naturals}{{\mathbb N}}

\begin{document}
\title{The Reznichenko property and the Pytkeev property in
hyperspaces}
\author{Ljubi\v{s}a D.R. Ko\v{c}inac}
\address{Department of Mathematics\\
Faculty of Sciences and Mathematics\\
University of Ni\v{s}\\
Vi\v{s}egradska 33\\
18000 Ni\v{s}, Serbia}
\email{lkocinac@ptt.yu}
\thanks{Supported by the Ministry of Science, Technology and Development, 
Republic of Serbia, grant N$^\text{o}$ 1233.}
\begin{abstract}
We investigate two closure-type properties, the Reznichenko property and
the Pytkeev property, in hyperspace topologies.
\end{abstract}
\date{}
\maketitle

\section{Introduction}

Let $X$ be a space (we suppose that all spaces are Hausdorff). For a
subset $A$ of $X$ and a family $\mathcal A$ of subsets of $X$ we put
$A^c=X\setminus A$ and $\mathcal A^c=\{A^c:A\in\mathcal A\}$. By $2^X$
we denote the family of all closed subsets of $X$. If $A$ is a subset of
$X$, then we write
\begin{center}
$A^{-} =\{F\in 2^X:F\cap A\neq\emptyset\}$,\\
$A^{+} =\{F\in 2^X:F\subset A\}$.
\end{center}

There are many known topologies on $2^X$. The most popular among them is
the Vietoris topology ${\sf V} = {\sf V^-} \vee {\sf V^+}$, where
the \emph{lower Vietoris topology} ${\sf V}^-$ is generated by all sets
$A^-$, $A\subset X$ open, and the \emph{upper Vietoris topology}
${\sf V}^+$ is generated by sets $B^+$, $B$ open in $X$. 

Let $\Delta$ be a subset of $2^X$. Then the \emph{upper $\Delta$-topology},
denoted by $\Delta^+$ (and first studied in abstract in \cite{poppe} and
then in \cite{dimaiohola}), is the topology whose subbase is the collection
\[
\{(D^c)^+:D \in \Delta\} \cup \{2^X\}.
\]
We consider only such subsets  $\Delta$ of $2^X$ which are {\bf closed
for finite unions} and {\bf contain all singletons}. In that case the
above collection is a base for $\Delta^+$ because we have
\[
(D_1^c)^+ \cap (D_2^c)^+= (D_1^c\cap D_2^c)^+= ((D_1\cup D_2)^c)^+ \, \,
\mbox{ and } \, \, D_1\cup D_2 \in \Delta.
\]
Two important special cases, in which we are especially interested in this
paper, are $\Delta= \mathbb F(X)$ -- the family of all finite subsets of
$X$, and $\Delta = \mathbb K(X)$ -- the collection of compact subsets of
$X$.  The $\mathbb F(X)^+$-topology will be denoted by ${\sf Z}^+$ and the
$\mathbb K(X)^+$-topology by ${\sf F^+}$. The ${\sf F}^+$-topology is
known as the \emph{upper Fell topology} (or the \emph{co-compact topology})
\cite{fell}. The \emph{Fell topology} ${\sf F}$ on $2^X$ is the
topology ${\sf V}^- \vee {\sf F^+}$. 

We investigate two closure type properties of $2^X$, the Reznichenko
property and the Pytkeev property, which are intermediate between
sequentiality and the countable tightness property that have been studied
in \cite{chv} and \cite{hou}.

Let us fix some terminology and notation that we need.

Let ${\mathcal A}$ and ${\mathcal B}$ be sets whose members are families
of subsets of an infinite set $X$. Then (see \cite{coc1}, \cite{coc2}):\\

\noindent $\sone({\mathcal A},{\mathcal B})$ denotes the selection
principle:
\begin{quote}
For each sequence $(A_n:n\in\naturals)$ of elements of ${\mathcal A}$
there is a sequence $(b_n:n\in\naturals)$ such that for each $n$
$b_n\in A_n$, and $\{b_n:n\in\naturals\}$ is an element of ${\mathcal B}$.
\end{quote}

\noindent $\sfin({\mathcal A},{\mathcal B})$ denotes the selection
hypothesis:
\begin{quote}
For each sequence $(A_n:n\in\naturals)$ of elements of ${\mathcal A}$
there is a sequence $(B_n:n\in\naturals)$ of finite sets such that
for each $n$ $B_n\subset A_n$, and $\bigcup_{n\in\naturals}B_n$ is
an element of ${\mathcal B}$.
\end{quote}

If $\mathcal C$ is a family of subsets of a space $X$ then an open cover
$\mathcal U$ of $X$ is called a \emph{$\mathcal C$-cover} if each $C\in
\mathcal C$ is contained in an element of $\mathcal U$.
$\mathbb F(X)$-covers and $\mathbb K(X)$-covers are customarily called
\emph{$\omega$-covers} and \emph{$k$-covers}, respectively.
We suppose that $\mathcal C$-covers we consider are non-trivial, i.e.
that $X$ does not belong to the cover.

\medskip
For a topological space $X$ and a point $x\in X$ we denote:
\begin{itemize}
\item[$1$.] $O\mathcal C$ -- the family of (open) $\mathcal C$-covers of
$X$.
\item[$2$.] $\Omega$ -- the family of $\omega$-covers of $X$.
\item[$3$.] $\mathcal K$ -- the family of $k$-covers of $X$.
\item[$4$.] $\Omega_x =\{A \subset X: x\in \overline{A}\setminus A\}$.
\end{itemize}
We also need the notion of groupability (see \cite{coc7}).

\smallskip
A countable $\mathcal C$-cover ${\mathcal U}$ of a space $X$ is
\emph{groupable} if there is a partition $({\mathcal U}_n:n\in
\naturals)$ of ${\mathcal U}$ into pairwise disjoint finite sets such
that for each $C\in\mathcal C$, for all but finitely many $n$ there is
a $U\in{\mathcal U}_n$ such that $C\subset U$.
A countable set $A\in\Omega_x$ is \emph{groupable} if there is a partition
$(B_n:n\in\naturals)$ of $A$ into finite sets such that each neighborhood
of $x$ has nonempty intersection with all but finitely many $B_n$.

\noindent Let us denote:
\begin{itemize}
\item[$5$.] $O\mathcal C^{gp}$ -- the family of groupable
$\mathcal C$-covers of $X$.
\item[$6$.] $\Omega^{gp}$ -- the family of groupable $\omega$-covers of
$X$.
\item[$7$.] $\mathcal K^{gp}$ -- the family of groupable $k$-covers of
$X$.
\item[$8$.] $\Omega_x^{gp}$ -- the family of groupable elements of
$\Omega_x$.
\end{itemize}
Let us mention that considering the space $2^X$, $\Delta \subset 2^X$ and
$S\in 2^X$ we shall use the symbol $\Delta^+_S$ to denote $\Omega_S$
with respect to $\Delta^+$-topology on $2^X$. We use
$\Omega_S$ and $\mathcal K_S$ studying the ${\sf Z}^+$- and
${\sf F}^+$-topology on $2^X$. Similarly for $(\Delta^+_S)^{gp}$
$\Omega_S^{gp}$ and $\mathcal K_S^{gp}$.

\section{ The Reznichenko property}

In 1996, Reznichenko introduced (in a seminar at Moscow State University)
the following property for a space $X$: For each $A\subset X$ and each
$x\in \overline{A}\setminus A$ there is a countable infinite family
$\mathcal A$ of finite pairwise disjoint subsets of $A$ such that each
neighborhood of $x$ meets all but finitely many elements of $\mathcal A$.
It is the same as to say that each countable element of $\Omega_x$ is a
member of $\Omega^{gp}_x$. This property was studied further in
\cite{malyhintironi} (under the name \emph{weakly Fr\'echet-Urysohn}),
\cite{fedeliledonne}. In \cite{ks} and \cite{coc7} this property was
called the \emph{Reznichenko property} and function spaces $\cpx$ having
this property were studied (see also \cite{sakai2}). In \cite{agt} this
property was considered in function spaces $\ckx$.

Evidently, if a space $X$ has the Reznichenko property then it has
countable tightness.

\medskip
Here we investigate the Reznichenko property in hyperspaces.

\begin{theorem}\label{rezdelta} For a space $X$ and a family $\Delta
\subset 2^X$ (closed for finite unions and containing all singletons)
the following statements are equivalent:
\begin{itemize}
\item[$(1)$] $(2^X,\Delta^+)$ has the Reznichenko property;
\item[$(2)$] For each open set $Y \subset X$ and each open $\Delta$-cover
$\mathcal U$ of $Y$ there is a sequence $(\mathcal U_n:n\in\naturals)$ of
finite pairwise disjoint subsets of $\mathcal U$ such that each $D\in
\Delta$ belongs to some $U\in\mathcal U_n$ for all but finitely many $n$.
\end{itemize}
\end{theorem}
$\pf$ $(1) \Rightarrow (2)$: Let $Y$ be an open subset of $X$ and let
$\mathcal U$ be an open $\Delta$-cover of $Y$. Then $\mathcal A:=
\mathcal U^c$ is a subset of $2^X$ and $Y^c\in Cl_{{\Delta}^+}(\mathcal
A)$. Indeed, let $D\in \Delta$ be a subset of $Y$. There is a $U\in
\mathcal U$ such that $D\subset U \subset Y$ and thus $Y^c\subset U^c
\subset D^c$, i.e. $U^c\in (D^c)^+$ and $Y^c\in (D^c)^+$. So, $Y^c\in
(D^c)^+\cap \mathcal A$, that is  $Y^c\in Cl_{{\Delta}^+}(\mathcal A)$.
Apply (1) to find a sequence $(\mathcal A_n:n\in \naturals)$ of finite
pairwise disjoint subsets of $\mathcal A$ such that each
$\Delta^+$-neighborhood of $Y^c$ intersects $\mathcal A_n$ for all but
finitely many $n$. For each $n$ let $\mathcal U_n =\mathcal A_n^c$. The
sets $\mathcal U_n\subset \mathcal U$ are pairwise disjoint (because
$\mathcal A_n$'s are) and witness for $\mathcal U$ that $Y$ satisfies
(2). Let $D\subset Y$ be an element from $\Delta$. Then $(D^c)^+$ is a
$\Delta^+$-neighborhood of $Y^c$, so that there is $n_0$ such that
$(D^c)^+\cap \mathcal A_n\neq\emptyset$ for each $n>n_0$. So, for each
$n>n_0$ there exists a set $A_n\in\mathcal A_n$ with $A_n\subset D^c$,
i.e. $D \subset A_n^c\in\mathcal U_n$. This means that (2) holds.\\

$(2) \Rightarrow (1)$: Let $\mathcal A$ be a subset of $2^X$ and
$S\in 2^X$ a point such that $S\in Cl_{\Delta^+}(\mathcal A) \setminus
\mathcal A$. Then $\mathcal U :=\mathcal A^c$ is a (non-trivial)
$\Delta$-cover of the open set $S^c\subset X$. Apply (2) to $S^c$ and
$\mathcal U$. One can choose a sequence $(\mathcal V_n:n\in\naturals)$
of finite pairwise disjoint subsets of $\mathcal U$ such that for each
$D\subset Y$ belonging to $\Delta$ for all but finitely many $n$ there is
a $V\in\mathcal V_n$ with $D\subset V$. Let for each $n$, $\mathcal B_n =
\mathcal V_n^c$. Then the collection $\{\mathcal B_n:n\in\naturals\}$
witnesses (1). Let $(D^c)^+$ be a $\Delta^+$-neighborhood of $S$. Then
$S\subset D^c$ implies $D\subset S^c$ so that there is $m\in\naturals$
such that for all $n>m$ there exists a member $B_n\in\mathcal B_n$ with
$D\subset B_n$, and consequently $B_n^c \in (D^c)^+$. This shows that
$(D^c)^+\cap \mathcal B_n\neq\emptyset$ for all but finitely many $n$
and completes the proof of the theorem.
$\epf$                                

\bigskip
In fact, we shall prove a general result regarding the Reznichenko
property in hyperspaces.

\begin{theorem}\label{rezgeneral1} Let $X$ be a space and let $\Delta$ and
$\Sigma$ be subsets of $2^X$ closed for finite unions and containing all
singletons. Then the following statements are equivalent:
\begin{itemize}
\item[$(1)$] $2^X$ satisfies $\sone(\Delta^+_A,(\Sigma^+_A)^{gp})$ for each
$A\in 2^X$;
\item[$(2)$] Each open set $Y \subset X$ satisfies $\sone(O\Delta,
O\Sigma^{gp})$.
\end{itemize}
\end{theorem}
$\pf$ $(1) \Rightarrow (2)$: Let $(\mathcal U_n:n\in\naturals)$ be a
sequence of $\Delta$-covers of $Y$. Then $(\mathcal U_n^c:n\in\naturals)$
is a sequence of subsets of $2^X$ and $Y^c\in Cl_{{\Delta}^+}(\mathcal
U_n^c)$ for each $n\in\naturals$. Applying (1) we find a sequence
$(U_n^c:n\in\naturals)$ such that for each $n$ $U_n\in\mathcal U_n$ and
$\mathcal F= \{U_n^c:n\in\naturals\}\in (\Sigma^+_{Y^c})^{gp}$. There is a
partition $\mathcal F=\bigcup_{n\in\naturals}\mathcal F_n$ of $\mathcal
F$ such that each ${\Sigma}^+$-neighborhood of $Y^c$ meets all but
finitely many sets $\mathcal F_n$. For each $n$ let $\mathcal V_n =
\mathcal F_n^c$. The sets $\mathcal V_n$ are pairwise disjoint, finite
subsets of $\{U_n:n\in\naturals\}$ and show that this set is a groupable
$\Sigma$-cover of $Y$. Let $S\subset Y$ be
a member of $\Sigma$. Then $(S^c)^+$ is a $\Sigma^+$-neighborhood of
$Y^c$, so that there is $n_0$ such that $(S^c)^+\cap \mathcal F_n\neq
\emptyset$ for each $n>n_0$. So, for each $n>n_0$ there exists a set
$F_n\in\mathcal F_n$ which is a subset of $S^c$, i.e. $S\subset F_n^c\in
\mathcal V_n$. This means that $\{U_n:n\in\naturals\}$ is a groupable
$\Sigma$-cover of $Y$, hence (2) holds.\\

$(2) \Rightarrow (1)$: Let $(\mathcal A_n:n\in\naturals)$ be a sequence
of subsets of $2^X$ such that a point $E\in 2^X$ belongs to 
$Cl_{\Delta^+}(\mathcal A_n)\setminus \mathcal A_n$ for each $n$. For
each $n$ put $\mathcal U_n =\mathcal A_n^c$. Apply (2) to the open set
$E^c$ and the sequence $(\mathcal U_n:n\in\naturals)$ of $\Delta$-covers
of $E^c$. We choose a sequence $(U_n:n\in\mathcal U_n)$ such that for
each $n$ $U_n\in\mathcal U_n$ and $\mathcal G= \{U_n:n\in\naturals\}$ is
a groupable $\Sigma$-cover of $E^c$. Suppose that the partition $\mathcal
G=\bigcup_{n\in\naturals}\mathcal G_n$ witnesses groupability of
$\mathcal G$. Let for each $n$, $A_n=U_n^c\in \mathcal A_n$ and $\mathcal
B_n= \mathcal G_n^c$. Then the collection $\{\mathcal B_n:n\in\naturals
\}$ witnesses that $\{A_n:n\in\naturals\}$ a groupable element of
$\Sigma^+_E$. Indeed, let $(S^c)^+$
be a $\Sigma^+$-neighborhood of $E$. Then $E\subset S^c$ implies $S
\subset E^c$ so that there is $k$ such that for all $n>k$ there exists a
member $G_n\in\mathcal G_n$ with $S\subset G_n$, and consequently $G_n^c
\in (S^c)^+$. This shows that $(S^c)^+\cap \mathcal B_n\neq\emptyset$
for all but finitely many $n$. 
$\epf$                                

\bigskip
The condition (2) in this theorem can be called the \emph{$(\Delta^+,
\Sigma^+)$-selectively Reznichenko property} (of $2^X$).

\bigskip
Let us recall that a space $X$ is said to have \emph{countable strong
fan tightness} \cite{sakai1} if for each sequence $(A_n:n\in\naturals)$
of subsets of $X$ and each $x\in \bigcap_{n \in \naturals}\overline{A}_n$
there is a sequence $(x_n:n\in \naturals)$ such that for each $n$
$x_n\in A_n$ and $x \in \overline{\{x_n:n\in \naturals\}}$, i.e. if for
each $x\in X$ the selection hypothesis $\sone(\Omega_x,\Omega_x)$ is
satisfied.

It was shown in \cite{dimkocmec} that $\zplusx$ (resp. $\fplusx$) has
countable strong fan tightness if and only if each open
set $Y\subset X$ satisfies $\sone(\Omega,\Omega)$ (resp. $\sone(\mathcal
K,\mathcal K)$).

\smallskip
According to  \cite{coc7} a space $X$ satisfies $\sone(\Omega,\Omega^{gp})$
if and only if all finite powers of $X$ have the \emph{Gerlits-Nagy
property} \cite{gerlitsnagy} (equivalently, the Hurewicz property as
well as the Rothberger property).
Recall that $X$ has the \emph{Rothberger property} \cite{rothberger38}
if it
satisfies $\sone(\mathcal O,\mathcal O)$, where $\mathcal O$ is the family
of open covers of $X$. $X$ has the \emph{Hurewicz property} \cite{hurewicz} if
for each sequence $(\mathcal U_n:n\in\naturals)$ of open covers of $X$
there is a sequence $(\mathcal V_n:n\in\naturals)$ such that for each
$n$ \, $\mathcal V_n$ is a finite subset of $\mathcal U_n$ and each $x
\in X$ belongs to $\cup\mathcal V_n$ for all but finitely many $n$.

Letting in Theorem \ref{rezgeneral1} $\Delta=\Sigma= \mathbb F(X)$ (resp.
$\Delta=\Sigma= \mathbb K(X)$) we obtain the following two important
corollaries.

\begin{corollary}\label{rezzet1} For a space $X$ the following statements
are equivalent:
\begin{itemize}
\item[$(1)$] $\zplusx$ has the Reznichenko property and countable
strong fan tightness; 
\item[$(2)$] If $Y$ is an open subset of $X$, then all finite powers of
$Y$ have the Gerlits-Nagy property.
\end{itemize}
\end{corollary}

\begin{corollary}\label{rezfell1} For a space $X$ the following are
equivalent:
\begin{itemize}
\item[$(1)$] $\fplusx$ has the Reznichenko property and countable strong
fan tightness; 
\item[$(2)$] Each open set $Y \subset X$ satisfies $\sone(\mathcal K,
\mathcal K^{gp})$.
\end{itemize}
\end{corollary}

\bigskip
Similarly to the proof of Theorem \ref{rezgeneral1} we can prove

\begin{theorem}\label{rezgeneralfin} Let $X$ be a space and let $\Delta$
and $\Sigma$ be subsets of $2^X$ closed for finite unions and containing
all singletons. Then the following statements are equivalent:
\begin{itemize}
\item[$(1)$] $2^X$ satisfies $\sfin(\Delta_A^+,(\Sigma^+_A)^{gp})$ for each
$A\in 2^X$;
\item[$(2)$] Each open set $Y \subset X$ satisfies $\sfin(O\Delta,
O\Sigma^{gp})$.
\end{itemize}
\end{theorem}

A space $X$ is said to have \emph{countable fan tightness} \cite{arh1},
\cite{arhbook} if for each
$x\in X$ it satisfies $\sfin(\Omega_x,\Omega_x)$, i.e. if whenever $(A_n:n\in
\naturals)$ is a sequence of subsets of $X$ and $x\in\bigcap_{n\in\naturals}
\overline{A}_n$ there are finite sets $B_n\subset A_n$, $n\in\naturals$,
such that $x \in \overline{\bigcup_{n\in\naturals}B_n}$.

In \cite{dimkocmec} it was proved that $\zplusx$ (resp. $\fplusx$) has
countable fan tightness if and only if each open set $Y\subset X$
satisfies $\sfin(\Omega,\Omega)$ (resp. $\sfin(\mathcal K,\mathcal K)$.
A result in \cite{coc7} states that all finite powers of of a space $X$
have the Hurewicz property if and only if $X$ belongs to the class
$\sfin(\Omega,\Omega)$. 

As the corollaries of Theorem \ref{rezgeneralfin} we have:

\begin{corollary}\label{rezzetfin} For a space $X$ the following
statements are equivalent:
\begin{itemize}
\item[$(1)$] $\zplusx$ has both the Reznichenko property and countable
fan tightness; 
\item[$(2)$] If $Y$ is an open subset of $X$, then each finite power of
$Y$ has the Hurewicz property.
\end{itemize}
\end{corollary}

\begin{corollary}\label{rezfellfin} For a space $X$ the following are
equivalent:
\begin{itemize}
\item[$(1)$] $\fplusx$ has the Reznichenko property and countable fan
tightness;
\item[$(2)$] Each open set $Y \subset X$ satisfies $\sfin(\mathcal K,
\mathcal K^{gp})$.
\end{itemize}
\end{corollary}            

\section{The Pytkeev property}

A space $X$ has the \emph{Pytkeev property} if for each $A\subset X$ and
each $x\in\overline{A\setminus \{x\}}$ there is a countable collection
$\{A_n:n\in\naturals\}$ of (countable) infinite subsets of $A$ which is
a $\pi$-network at $x$, i.e. each neighborhood of $x$ contains some
$A_n$. This property was introduced in \cite{pytkeev} and then studied
in \cite{malyhintironi} (where the name Pytkeev space was used) and
\cite{fedeliledonne}. Pytkeev's property in function spaces
$\cpx$ was studied in \cite{sakai2}.

Every (sub)sequential space has the Pytkeev property \cite[Lemma 2]
{pytkeev} and every Pytkeev space has the Reznichenko property
\cite[Corollary 1.2]{malyhintironi}.

In this section we consider the Pytkeev property in hyperspaces with
$\Delta^+$-topologies. We begin by the following  general result.

\begin{theorem}\label{pytkeevgeneral} If $X$ is a space and $\Delta$ and
$\Sigma$ subsets of $2^X$ containing all singletons and closed for finite
unions, then the following are equivalent:
\begin{itemize}
\item[$(1)$] $2^X$ has the $(\Delta^+,\Sigma^+)$-Pytkeev property, i.e.
for each $\mathcal A\subset 2^X$ and each $S \in Cl_{\Delta^+}(\mathcal
A \setminus \{S\})$ there is a family $\{\mathcal B_n:n\in\naturals\}$ of
countable infinite subsets of $\mathcal A$ which is a $\pi$-network at
$S$ with respect to the $\Sigma^+$-topology. 
\item[$(2)$] For each open set $Y \subset X$ and each $\Delta$-cover
$\mathcal U$ of $Y$ there is a sequence $(\mathcal U_n:n\in\naturals)$
of infinite subsets of $\mathcal U$ such that $\{\cap\, \mathcal U_n:n
\in\naturals\}$ is a (not necessarily open) $\Sigma$-cover of $Y$.
\end{itemize}
\end{theorem}
$\pf$ $(1) \Rightarrow (2)$: Let $\mathcal U$ be a $\Delta$-cover of $Y$.
Then $\mathcal A :=\{U^c:U\in\mathcal U\}$ is a subset of $2^X$ and
$Y^c\in Cl_{\Delta^+}(\mathcal A)$. Apply (1) to find a sequence
$(\mathcal B_n:n\in\naturals)$ of countable infinite subsets of
$\mathcal A$ such that the set $\{\mathcal B_n:n\in\naturals\}$ ia a
$\Sigma^+$-$\pi$-network at $Y^c$. For each $n$ denote by $\mathcal U_n$
the subset $\{U:U^c\in\mathcal B_n \}$ of $\mathcal U$. We claim that
$\{\cap\, \mathcal U_n:n\in\naturals\}$ is a $\Sigma$-cover of $Y$. 
For any $S\subset Y$ with $S\in\Sigma$, the set $(S^c)^+$ is a
$\Sigma^+$-neighborhood of $Y^c$ and therefore there is $m\in\naturals$
such that $\mathcal B_m\subset (S^c)^+$. Therefore, for each $U^c\in
\mathcal B_m$ we have $U^c\subset S^c$, i.e. $S\subset \cap\{U:U^c\in
\mathcal B_m\}$. This means that the collection $\{\cap\,\mathcal U_n:n
\in\naturals\}$ is indeed a $\Sigma$-cover of $Y$.\\

$(2)\Rightarrow (1)$: Let $\mathcal A$ be a subset of $2^X$ and let
$S\in Cl_{{\Delta}^+}(\mathcal A \setminus \{S\})$. Then $S^c$ is an open
subset of $X$ and the family $\mathcal U \equiv \mathcal A^c:=\{A^c:A\in
\mathcal A\}$ is an open $\Delta$-cover of $S^c$. Apply (2) to the set
$S^c$ and its $\Delta$-cover $\mathcal U$ and choose infinite sets
$\mathcal U_n\subset \mathcal U$, $n\in\naturals$, such that the family
$\{\cap\, \mathcal U_n:n\in\naturals\}$ is a $\Sigma$-cover of $S^c$.
Put for each $n$
\[
\mathcal A_n=\{U^c:U\in\mathcal U_n\}.
\]
Then each $\mathcal A_n$ is an infinite subset of $\mathcal A$. We prove
that the collection $\{\mathcal A_n:n\in\naturals\}$ is a $\pi$-network
at $S\in (2^X,\Sigma^+)$.

Let $(E^c)^+$, $E\subset X$ and $E\in\Sigma$, be a ${\Sigma}^+$-neighborhood
of $S$. Then $S\subset E^c$ implies $E\subset S^c$ so that there is some
$i\in\naturals$ with $E\subset \cap\, \mathcal U_i$. Further we have
$\mathcal A_i\subset (E^c)^+$, i.e. (1) holds.
$\epf$

\bigskip
As consequences of this theorem we obtain the following two results.

\begin{corollary}\label{pytkeevzet} For a space $X$ the following
statements are equivalent:
\begin{itemize}
\item[$(1)$] $(2^X,{\sf Z}^+)$ has the Pytkeev property;
\item[$(2)$] For each open set $Y \subset X$ and each $\omega$-cover
$\mathcal U$ of $Y$ there is a sequence $(\mathcal U_n:n\in\naturals)$
of infinite countable subsets of $\mathcal U$ such that $\{\cap\,
\mathcal U_n:n\in\naturals\}$ is an (not necessarily open) $\omega$-cover
of $Y$.
\end{itemize}
\end{corollary}

\begin{corollary}\label{pytkeevfell} For a space $X$ the following
statements are equivalent:
\begin{itemize}
\item[$(1)$] $(2^X,{\sf F}^+)$ has the Pytkeev property;
\item[$(2)$] For each open set $Y \subset X$ and each $k$-cover $\mathcal
U$ of $Y$ there is a sequence $(\mathcal U_n:n\in\naturals)$ of countable
infinite subsets of $\mathcal U$ such that $\{\cap\, \mathcal U_n:n\in
\naturals\}$ is a (not necessarily open) $k$-cover of $Y$.
\end{itemize}
\end{corollary}

\begin{problem}
If $\fplusx$ has the Pytkeev property, is $\fplusx$ sequential? 
What about $\zplusx$?
\end{problem}

\noindent{\bf Remark A} It is known that the space $(2^X,{\sf F})$ is
compact with no assumptions on $X$ and that a Hausdorff space $X$ is
locally compact if and only if $(2^X,{\sf F})$ is Hausdorff. According
to a result from \cite{malyhintironi}, each compact Hausdorff space of
countable tightness is a Pytkeev space. On the other hand, in \cite{hou}
(see also \cite{chv}) it was shown that for a locally compact Hausdorff
space $X$ the tightness of $(2^X,{\sf F})$ is countable if and only if
$X$ is hereditarily separable and hereditarily Lindel\"of. So we have:

\begin{theorem} For a locally compact Hausdorff space $X$ the following
assertions are equivalent:
\begin{itemize}
\item[$(a)$] $(2^X,{\sf F})$ has countable tightness;
\item[$(b)$] $(2^X,{\sf F})$ has the Reznichenko property;
\item[$(c)$] $(2^X,{\sf F})$ has the Pytkeev property;
\item[$(d)$] $X$ is both hereditarily separable and hereditarily
Lindel\"of.
\end{itemize}
\end{theorem}

\noindent{\bf Remark B} According to results from \cite{balogh} and
\cite{dow} one can suppose that in some models of ${\sf ZFC}$ (in which
each compact Hausdorff space of countable tightness is sequential) we
have: For a locally compact Hausdorff space $X$ each of the conditions
(a)--(d) in the previous theorem is equivalent to the assertion $(2^X,
{\sf F})$ is a sequential space.

\bigskip
\noindent{\bf Remark C} Let us consider a selective version of the
Pytkeev property. Call a space $X$ \emph{selectively Pytkeev} if for
each sequence $(A_n:n\in\naturals)$ of subsets of $X$ and each point
$x\in \bigcap_{n\in\naturals}\overline{A_n\setminus \{x\}}$
there  is an infinite family $\{B_n:n\in\naturals\}$ of countable
infinite sets which is a $\pi$-network at $x$ and such that for each $n$
$B_n\subset A_n$. Then one can prove the statements similar to those in
Theorem \ref{pytkeevgeneral} and Corollaries \ref{pytkeevzet} and
\ref{pytkeevfell}.

\end{document}